\documentclass[a4paper,11pt]{article}
\usepackage{amsmath,amscd,amsfonts,amssymb}
\newtheorem{theorem}[equation]{Theorem}

\newtheorem{question}[equation]{Question}

\newtheorem{corollary}[equation]{Corollary}
\newtheorem{remark}[equation]{Remark}
\newtheorem{definition}[equation]{Definition}

\def\C{\mathbb{C}}
\def\Q{\mathbb{Q}}

\def\Zc{{\cal Z}}

\def\N{\mathbb{N}}

\def\P{\mathbb{P}}
\def\qed{\hfill$\Box$\s}
\def\s{\vskip10pt}

\date{May 2015}

\begin{document}

\title{\bf Leray spectral sequence for complements of certain arrangements of smooth submanifolds }
\author{
Andrzej Weber\thanks{Supported by NCN grant 2013/08/A/ST1/00804
}
\\
\small Department of Mathematics of Warsaw University\\
\small Banacha 2, 02-097 Warszawa, Poland\\
\small aweber@mimuw.edu.pl}

\maketitle

Spectral sequence argument quite often appears in study the hyperplane arrangements or configuration spaces (see for example \cite{BE, DJL, DSY, DS, Je, JOS, To}). Very often the spectral sequence can be identified with the Leray spectral sequence of the embedding $j:U\hookrightarrow X$ of an open subset to a compact complex variety
\begin{equation} E^{p,q}_2=H^p(X;R^qj_*{\cal L})\;\Rightarrow\; H^{p+q}(U;{\cal L})\,.\label{leray}\end{equation}
 When the coefficient system has geometric nature, then by Deligne the entries of the spectral sequence are equipped with the weight filtration
 in the sense of the mixed Hodge structure \cite{De1, De2} or in the sense of the Weil sheaves \cite{De0, BBD}.
 For the constant coefficients $\Q_U$ it happens quite often that the entries of $E_2$ table are of pure weight. That is so for:
 \begin{description}
 \item[(1)] complements of the classical hyperplane arrangements in $\P^n$,
 \item[(2)] complements of elliptic hyperplane arrangements,
 \item[(3)] in general: complements of hyperplane-like arrangements.
 \end{description}
 Purity implies the degeneration of the spectral sequence on $E_3$. The argument of purity can be applied in a much more general setup. We can add to the list above
 \begin{description}
 \item [(4)]configuration space of $n$ points in a projective manifold $Y$ considered as an open subset of $X=Y^n$
\end{description}
 It was already noticed by Totaro \cite{To} that if $X=Y^n$ for a smooth compact algebraic variety $Y$ and $U$, the set of tuples consisting of pairwise distinct points, i.e. when $U$ is the configuration space of $n$ points in $Y$, then the Leray spectral sequence has only one nontrivial differential.

 Our goal is to add to the list of spaces with pure sheaves $R^qj_*\Q_U$ another class of ,,arrangements'', which we call ,,admissible arrangements of submanifolds''. The members of the arrangements are smooth subvarieties of $X$, but we allow some kind of singular intersections. For example we allow that two subvarieties are tangent along a smooth subvariety of certain codimension, see Definition \ref{defi}. We add to the list of open sets for which Leray sequence degenerates another item\begin{description}
 \item [(5)] complements of admissible arrangements of submanifolds in a compact smooth algebraic variety.
 \end{description}
 The concept of weight is central for this paper. Talking about weight filtration in the cohomology of complex algebraic variety one obviously uses Deligne construction \cite{De1, De2}. While dealing with sheaves there are (at least) two approaches possible. It is natural to work in the category of mixed Hodge modules developed by M.~Saito \cite{Sa}. Nevertheless we feel more comfortable using the reduction to $\Q_\ell$-sheaves on a variety over a finite characteristic field, as in \cite{De0, BBD}. If we say that a constructible sheaf of $\Q$-vector spaces is pure of certain weight we mean that after tensoring with $\Q_\ell$ and after a good reduction to a finite characteristic $p$ with $(p,\ell)=1$, the sheaf is pure in the sense of \cite[\S 5.1]{BBD} (i.e.~in the sense of the action of the Frobenius automorphism on stalks and costalks).

Many results about arrangements of hyperplanes can be obtained purely topologically, one can say ,,elementarily''.
Sometimes the weight argument is hidden somewhere in the references.
Our goal is to review some basic mechanism and extend a bit the situation when it can be applied. The weight argument works for cohomology with rational coefficients. Some of the results might hold for integer homology, but the weight argument is too weak to handle torsion.

For the twisted coefficients, i.e.~when the local system ${\cal L}$ is not constant, the Leray spectral sequence allows to generalize proofs done in particular situations. The conditions on the monodromy eigenvalues {\bf(Mon)} of \S\ref{monodromia} implies that the stalks $R^kj_*{\cal L}$ vanish at points $x\in X\setminus U$. Therefore
$$H^*(U;{\cal L})=H^*_c(U;{\cal L})\,,$$
where $H_c^*(-)$ denotes cohomology with compact supports. To prove that no use of Hodge theory is needed. The situation when the stalks of $R^kj_*{\cal L}$ do not vanish is more interesting. Another condition, called $\bf (Mon)^*$ in the classical setup, guaranties that the local cohomology is generated by logarithmic forms. This condition was studied e.g.~in \cite[Theorem 4.1]{STV}. If the local system $\cal L$ is equipped with a pure Hodge structure it seems that the stalks should be pure. The problem lies in the choice of a convenient category in which purity is understood. A natural candidate would be the category of {\it mixed Hodge modules} of M. Saito \cite{Sa}. Making this choice we would have to add highly nontrivial technical introduction. We limit ourself just to few formal remarks in the final section. The supporting evidence is that in the case of normal crossing divisor and unitary local system the Leray spectral sequence degenerates. This was proven by Timmerscheidt in \cite{Ti}.
\s
This note was inspired by the Workshop on Configuration Spaces in Cortona. I'm grateful to Alex Suciu for conversations and encouragement.
\section{The results}
We will consider the following situation. Let $X$ be a smooth compact algebraic variety.
Let $r$ be a natural number. We
define an admissible arrangement of codimension $c$ submanifolds inductively.

\begin{definition}\label{defi}An arrangement of distinct submanifolds $\Zc=\{Z_i\}_{i=1,2,\dots s}$ is admissible of codimension $c$ in $X$ if
\begin{enumerate}
 \item each $Z_i$ is smooth of codimension $c$
 \item for any set of indices $I\subset\{1,2,\dots,s\}$ and $j\not\in I$ the set-theoretic intersection $$\left(\bigcup_{i\in I}Z_i\right)\cap Z_j$$ is empty or it is the union of an admissible codimension $c$ arrangement in each $Z_j$
\end{enumerate}\end{definition}
The definition of admissible arrangement of submanifolds is constructed in the way which allows to apply the argument of ,,deletion-restriction'' of \cite[Proposition 2.1]{DS} or \cite[Lemma 2.1]{DJL}.
Clearly each arrangement of submanifolds satisfying the following condition is admissible:\begin{description}\it
\item[$\bigstar$] for any set of indices $I\subset\{1,2,\dots,s\}$ the set-theoretic intersection $\bigcap_{i\in I}Z_i$ is empty or each irreducible component of $\bigcap_{i\in I}Z_i$ is a smooth subvariety
of codimension $c\ell$ for some $\ell\in \N$.
\end{description}

To prove that $\bigstar$ implies the condition of Definition \ref{defi} assume that $Z_i$'s are connected and they do not repeat.
Then the intersection $\left(\bigcup_{i\in I}Z_i\right)\cap Z_j=\bigcup_{i\in I}(Z_i\cap Z_j)\subset Z_j$ is the sum of smooth irreducible components. The members of the arrangement $\{Z_i\cap Z_j\}_{i\in I}$, if not empty, have codimension $\leq 2c$ in $X$ and by the assumption $\bigstar$ the codimension is a multiple of $c$. Therefore it has to be equal to $2c$. Hence $Z_i\cap Z_j$ has codimension
$c$ in $Z_j$ or it is empty. Some members of this arrangement may repeat, but we replace the arrangement by the set of irreducible components of intersections $Z_i\cap Z_j$. It again sa\-tis\-fies the condition $\bigstar$, therefore it satisfies the inductive step of Definition \ref{defi}.
\begin{remark}\rm A condition related to $\bigstar$ is considered in \cite[\S3.1]{DSY}: the intersection
$\bigcap_{i\in I}Z_i$ is assumed to be smooth. Then the topological components coincide with the irreducible components. On the other hand there is no assumption on the dimensions of intersections. The paper \cite{DSY} is strictly related to ours.
\end{remark}

In particular when $c=1$ then $\bigstar$ reduces to the condition
\begin{description}\item[$\bullet$]\it $\{Z_i\}$ is an arrangement of smooth hypersurfaces $Z_i$, such that for each set of indices $I$ any irreducible component of the intersection $\cap_{i\in I}Z_i$ is smooth.
\end{description}

Note that in \cite[\S3]{STV} there were considered hyperplane-like arrangements, i.e. arrangement of hypersurfaces which locally look like hyperplane arrangements. The hyperplane-like arrangements satisfy the Definition \ref{defi} with $c=1$. Our definition allow are more complicated examples: in the normal direction the intersection may look like the arrangement of three curves $x=0$, $x=y^2$ and $y=0$ on the plane. Another example of an admissible arrangement of hyperplanes is the sum of two surfaces in $\C^3$ given by the equations $z=0$ and $z=xy$. Again this is a picture in the normal slice, and the arrangement itself may be of high dimension.

\begin{remark}\rm Combinatorial nature of Leray spectral sequence for admissible arrangements of submanifolds has lead the authors of \cite{DSY} to a definition of the combinatorial cover and a related spectral sequence. The relation of their construction with Leray spectral sequence is explained in [{\it loc.cit}, \S2.2]. In [{\it loc.cit}, Theorem 3.1] a cover indexed by components of intersections $\bigcap_{i\in I}Z_i$ is constructed. The authors take sufficiently small tubular neighbourhoods of these components. This particular covering gives rise to a spectral sequence which is isomorphic to ours after some re-indexing of entries.\end{remark}

We will prove a theorem about degeneration of the spectral sequence for cohomology with trivial coefficients:

\begin{theorem}\label{t1} Let $\Zc=\{Z_i\}_{i=1,2,\dots s}$ be an admissible arrangement of codimension $c$ submanifolds in $X$. Then the Leray spectral sequence of the inclusion $$j:U=X\setminus \bigcup_{i\in\{1,2,\dots s\}} Z_i\longrightarrow X$$
$$E^{p,q}_2=H^p(X;R^qj_*\Q_U)\;\Rightarrow\; H^{p+q}(U;\Q)\,.$$
has all differentials $d_r$, $r\geq2$ vanishing except from $d_{2c}$.\end{theorem}
The result follows from the local computation:
\begin{theorem}\label{t2} With the notation as above the sheaves $R^kj_*\Q_U=0$ for $k$ not divisible by $2c-1$ and for $k={(2c-1)}\ell$ the sheaf $R^kj_*\Q_U$ is pure of weight $2c\ell=\frac{2c}{2c-1}k$.\end{theorem}
After further work one can show:
\begin{theorem}\label{t3} With the notation as above the sheaf $R^{(2c-1)\ell}j_*\Q_U$ is isomorphic to a direct sum of constant sheaves supported by smooth subvarieties of codimensions $c\ell$.\end{theorem}

We would like to stress that our degeneration result, Theorem \ref{t1}, specializes to many known cases mentioned in the introduction as the cases {\bf (1)--(4)}.

\section{Ideal situation: normal crossing divisor}

Let us start with recalling the Deligne
spectral sequence defining the weight filtration in the cohomology of the normal crossing divisor complement $X\setminus D$. The the sheaf of meromorphic forms on $X$ with logarithmic poles along $D$
admits the weight filtration defined by the number of factors in denominators. The iterated residue defines an isomorphism from the graded pieces of the filtration to the sheave of forms on the intersections of the divisors. One obtains a spectral sequence allowing to compute the cohomology of $X\setminus D$ from the cohomologies of the intersections. The main theorem says that the spectral sequence degenerates on $E_2$, hence the only thing one has to know is the first differential.

Let us be more precise. We fix the notation:
Let $X$ be a compact smooth complex algebraic variety. We assume
that $D\subset X$ is a smooth divisor with normal crossings
i.e.~the irreducible components of $D=\bigcup_{i=1}^m D_i$ are smooth and locally $D$ is
given by the equation $z_1z_2\dots z_k=0$ in a certain
coordinate system.
For a multindex $I=\{i_1,i_2,\dots,i_\ell\}$ define
$X_I=D_{i_1}\cap D_{i_2}\cap\dots\cap D_{i_\ell}$ and
$X^\ell=\coprod_{|I|=\ell}X_I$ with $X^0=X$.
Deligne considers the spectral sequence associated to the weight filtration in the complex of logarithmic forms \cite[(3.2.4.1)]{De1}. It converges to $H^*(X\setminus D)$.
The first table of the weight spectral sequence lies in the second quarter and has the form
$$_WE^{-p,q}_1=H^{q-2p}(X^p)\,,$$
\cite[(3.2.7)]{De1}.
Here is the example for the surface case:

\hfil $\begin{matrix}0&\to&H^0(X^2)&\to&H^2(X^1)&\to&H^4(X^0)&^{_4}\cr
 & & 0 &\to&H^1(X^1)&\to&H^3(X^0)&^{_3}\cr
 & & 0 &\to&H^0(X^1)&\to&H^2(X^0)&^{_2}\cr
 & & & & 0 &\to&H^1(X^0)&^{_1}\cr
 & & & & 0 &\to&H^0(X^0)&^{_0}\cr
 ^{_{-3}}& &^{_{-2}} & & ^{_{-1}} & & ^{_0} \cr\end{matrix}$
\s
\noindent
The differential $_Wd_1$ is the alternating sum of Gysin maps $$H^*(X_I)\to
H^{*+2}(X_{I\setminus\{i\}})$$ induced by the inclusions
$X_I\hookrightarrow X_{I\setminus \{i\}}$. Since we care about the Hodge
structure we should rather write
$H^*(X^I)(-1)\to
H^{*+2}(X^{I\setminus\{i\}})$, where $V(-1)$ denotes the Tate twist
 rising the weight by 2.
The differential $$_Wd_1:H^*(X^p)(-p)\to H^{*+2}(X^{p-1})(1-p) $$ might be nontrivial.
The weight principle says: the maps between Hodge structures of distinct weights vanish.
Therefore the second differential of
the spectral sequence is the zero map and the weight spectral sequence degenerates on $_WE_2$, \cite[(3.2.10)]{De1}.
The rows of the spectral sequence $_WE_2$ compute the graded pieces of the weight filtration in $H^*(U)$
\s
\def\ot{\leftarrow}
\noindent $Gr^W_{k+\ell} H^k(U)=$\s
 \hfill $=H_\ell\bigl(
 0
\ot\underbrace{ H^k(X^0)}_0
\ot\underbrace{ H^{k-2}(X^1)(-1)}_{1}
\ot\underbrace{ H^{k-4}(X^2)(-2)}_{2}
\ot\dots\;\bigr)\,.$

\noindent
As a matter of fact the weight spectral sequence after a reorganization of indices coincides with the
Leray spectral sequence of the inclusion
$j:X\setminus D \hookrightarrow X$.
We have $$R^qj_*\Q_U=\bigoplus_{|I|=q}\Q_{X_I}(-q)\,.$$
The Leray spectral sequence in this case has the following second table
$$E^{p,q}_2=H^p(X^q)(-q).$$
Here is the picture for the surface case
$$\begin{matrix}
 ^{_2} &H^0(X^2)(-2)&0&0&0&0\cr
 ^{_1} &H^0(X^1)(-1)&H^1(X^1)(-1)&H^2(X^1)(-1)&0&0\cr
 ^{_0} &H^0(X^0)&H^1(X^0)&H^2(X^0)&H^3(X^0)&H^4(X^0)\cr
 &^{_{0}} &^{_{1}} & ^{_{2}} & ^{3} &^{4}\cr\end{matrix}$$
The second differential
$$d_2:E^{p,q}_2=H^p(X^q)(-q)\to E^{p+2,q-1}_2=H^{p+2}(X^{q-1})(1-q)$$
preserves the weight: both entries have weight $p+2q$. Checking directly the definitions of both differentials we find, that $d_2={_Wd}_1$. According to the weight principle the Leray spectral sequence degenerates on $E_3$, \cite[(3.2.13)]{De1}.
\begin{remark}\rm Let $\cal L$ be an unitary local system on a normal crossing divisor complement. The degeneration of the Leray spectral sequence (\ref{leray})
 was proven in \cite{Ti}.\end{remark}

The goal of this paper is to show that a similar degeneration happens for an admissible arrangement of submanifolds when cohomology are taken with constant coefficients, i.e. Theorem \ref{t1}. In the final section we will briefly discuss a possible degeneration for nontrivial local systems, but we are bounded by the category in which our object live. We will assume that ${\cal L}$ is of geometric origin to be able to talk about purity.

\section{Proofs}

{\it Proof of Theorem \ref{t2}}. The first part of our proof is essentially the ,,deletion-restriction'' of \cite[Proposition 2.1]{DS} or \cite[Lemma 2.1]{DJL} rewritten in the sheaf language. We will prove a statement a bit stronger: the sheaves $R^kj_*\Q_U$ are pure and pointwise pure (compare the definition of purity \cite[5.1.8]{BBD}).
We proceed inductively with respect to the number of elements of the arrangement $\Zc$ and the dimension of $X$. If $\Zc=\{Z\}$, then there is a long exact sequence
$$\to H^{k-1}(U)\to H^{k-2c}(Z)(-c)\to H^k(X)\to H^k(U)\to\,,$$
 which localizes to an exact sequence of sheaves
$$\to R^{k-1}j_*\Q_U\to {\cal H}^{k-2c}(\Q_Z)(-c)\to {\cal H}^{k}(\Q_X)\to R^kj_*\Q_U\to\,,$$
Therefore
 $R^0j_*\Q_U=\Q_X$ and $R^{2c-1}j_*\Q_U=\Q_{Z}(-c)$.
The derived images vanish in the remaining degrees.
In particular the theorem holds when the dimension of $Z_i$ is zero. Assume now that the theorem holds for all families in varieties of dimension smaller then $n=\dim( X)$. We also assume that the theorem holds for families consisting of $s-1$ elements. Consider the Mayer-Vietoris exact sequence for the pair $U_1=X\setminus Z_1$, $U_2=X\setminus \bigcup_{i>1}Z_i$:
$$\to H^{k-1}(U)\to H^k(U_1\cup U_2)\to H^k(U_1)\oplus H^k(U_2)\to H^k(U)\to$$
which can be localized to the long exact sequence of sheaves
$$\to R^{k-1}j_*\Q_U\to R^kj'_*\Q_{U_1\cup U_2}\to R^kj_{1*}\Q_{U_1}\oplus R^kj_{2*}\Q_{U_2}\to R^kj_*\Q_U\to\,,$$
where $j_1,j_2,j'$ are inclusions of the sets $U_1$, $U_2$ and $U_1\cup U_2$ into $X$. Note that
$$R^kj'_*\Q_{U_1\cup U_2}=\left\{\begin{matrix}\Q_X&\text{for }k=0\,,\hfill\\
0&\text{for }0<k\leq 2c\,,\hfill\\
R^{k-2c}j''_*\Q_{Z_1\cap U_2}(-c)&\text{for }k>2c\,, \hfill\end{matrix}\right.$$
where $j'':Z_1\cap U_2\to X$ is the inclusion.
Topologically this formula can be justified by the following: locally, i.e. if $X$ was a ball, then $U_1\cup U_2$ would be the $2c$-fold suspension of $Z_1\cap U_2$. To compute the weight we notice that for $k>2c$, $x\in Y=Z_1\cap \bigcup_{i>1}Z_i$ by Alexander duality in $X$
$$\left(R^kj'_*\Q_{U_1\cup U_2}\right)_x={\cal H}_{2n-k-1}(Y)_x(n)\,,$$
where ${\cal H}_*$ is the sheaf of local homology\footnote{
Avoiding the sheaf language and forgetting Hodge structure we have $$\left(R^kj'_*\Q_{U_1\cup U_2}\right)_x\simeq
H^k(B_x\cap(U_1\cup U_2))\simeq H^{k}(\partial B_x\cap(U_1\cup U_2))\,,$$ where $B_x$ is a small ball around $x$. The stalk of the sheaf of local cohomology ${\cal H}_{\ell}(Y)_x$ is isomorphic to
$$H_{\ell}(B_x\cap Y, \partial B_x\cap Y)\simeq \bar H_{\ell-1}(\partial B_x\cap Y)\,.$$ By Alexander duality in $\partial B_x$ we have
$$\bar H^{k}(\partial B_x\cap(U_1\cup U_2))\simeq\bar H_{2n-k-2}(\partial B_x\cap Y)\simeq H_{2n-k-1}(B_x\cap Y, \partial B_x\cap Y)\,.$$}.
Again, by Alexander duality in $Z_1$, which is of dimension $n-c$,
$${\cal H}_{2n-k-1}(Y)_x(n)=\left(R^{k-2c}j''_*\Q_{Z_1\cap U_2}\right)_x(-c)\,.$$
Hence
$$\left(R^kj'_*\Q_{U_1\cup U_2}\right)_x\simeq \left(R^{k-2c}j''_*\Q_{Z_1\cap U_2}\right)_x(-c)\,.$$
for $k>2c$. This is a kind of suspension isomorphism.
Using the inductive hypothesis we find that $R^kj'_*\Q_{U_1\cup U_2}$ is nonzero only if $k=0$ or $k-1$ is divisible by $2c-1$ and the weight of $R^{(2c-1)\ell+1}j'_*\Q_{U_1\cup U_2}$ is equal to $2c\ell$. If $c>1$ the map
$R^kj'_*\Q_{U_1\cup U_2}\to R^kj_{1*}\Q_{U_1}\oplus R^kj_{2*}\Q_{U_2}$ vanishes by dimensional reasons. For $c=1$ we use the weight argument: the source is of weight $2k-1$, the target is of weight $2k$. We conclude that we obtain a short exact sequence of sheaves of weight $\frac{2c}{2c-1}k$,
$$0\to R^kj_{1*}\Q_{U_1}\oplus R^kj_{2*}\Q_{U_2}\to R^{k}j_*\Q_U\to R^{k+1}j'_*\Q_{U_1\cup U_2}\to 0\,.$$
(For $k=0$ the exact sequence is extended from the left by $R^0j'_*\Q_{U_1\cup U_2}\simeq \Q_X$, which is pure of weight 0.)
By inductive assumption the edge sheaves are pure and pointwise pure. Therefore the middle one is pure and pointwise pure.
For $k=0$ the exact sequence is extended from the left by $R^0j'_*\Q_{U_1\cup U_2}$, which is pure of weight 0.\qed

{\it Proof of Theorem \ref{t3}}. Again by induction we find that the sheaf $R^{(2c-1)\ell}j_*\Q_U$ fits to an exact sequence
$$0\to \bigoplus_\alpha \Q_{A_\alpha}(c\ell)\to R^{(2c-1)\ell}j_*\Q_U\to \bigoplus_\beta \Q_{B_\beta}(c\ell)\to 0\,.$$
Here $A_\alpha$ and $B_\beta$ are families of smooth varieties of codimensions $c\ell$.
By \cite[5.4.6]{BBD} $R^{(2c-1)\ell}j_*\Q_U$ decomposes into the sum of constant sheaves. (In the case of when the arrangement is homeomorphic to a vector space arrangement see \cite[III \S3.7]{GM} for a topological proof.) \qed

{\it Proof of Theorem \ref{t1}.} We follow the Totaro argument \cite[\S4]{To} (see also his remark at the beginning of p.1062) . If a sheaf $S$ is pure of weight $w$, the base space is compact, then the cohomology $H^p(X;S)$ is pure of weight $p+w$ (\cite[5.1.13]{BBD}). By Theorem \ref{t2} the entries of the Leray spectral sequence $E_2^{p,q}=H^p(X;R^qj_*\Q_U)$ are pure of weight $p+\frac{2c}{2c-1}q$ for $q$ divisible by $2c-1$. The remaining rows of the $E_2$ table are zero. The first possible differential is $$d_{2c}:E^{p,(2c-1)\ell}_{2}\longrightarrow E^{p+2c,(2c-1)\ell-2c+1}_2\,,$$
that is
$$d_{2c}:H^p(X;R^{(2c-1)\ell}j_*\Q_U)\longrightarrow H^{p+2c}(X;R^{(2c-1)(\ell-1)}j_*\Q_U)\,.$$
The weights of both entries are equal to $p+2c\ell=(p+2c)+2c(\ell-1)$ and the differential $d_{2c}$ can be nontrivial.
Further differentials hit the entries $E^{p+k,(2c-1)\ell-k+1}_{2c+1}$ which are subquotients of $H^{p+k}(X;R^{(2c-1)\ell-k+1}j_*\Q_U)$ and have weights $$(p+k)+\frac{2c}{2c-1}((2c-1)\ell-k+1)=p+2c\ell+\frac{2c-k}{2c-1}<p+2c\ell$$
for $k>2c$. Therefore the higher differentials vanish.\qed

\section{Twisted coefficients}\label{monodromia}
We would like to conclude with some remarks concerning the cohomology with coefficients in a local system. We restrict our attention to the arrangements of hypersurfaces ($c=1$), when the local monodromy might be nontrivial.
If ${\cal L}$ is a local system on $U=X\setminus D$ the Leray spectral sequence relates the local invariants of intersections of divisors with the global cohomology. Consider the case of a hyperplane-like arrangement, i.e.~we assume as in \cite{STV} that in some local analytic coordinates the hypersurafces are given by linear equation. Under the
condition that:\begin{itemize}\item\it $\bf(Mon)$:
 $\dim({\cal L})=1$ and for each $x\in D$ the product of the monodromies associated to the divisors passing through $x$ is not equal to one.
 \end{itemize}
 one obtains the vanishing of local cohomology $$(R^kj_*{\cal L})_x=H^k(U\cap B_{x} ;{\cal L})=0$$ for any $k$ at any point belonging to $D$, see \cite[Lemma 5.2]{DJL}. (Here $B_{x}$ is a sufficiently small ball centered in $x$.)
For abelian local systems of higher rank the generalization of {\bf(Mon)} would demand that the products of all possible eigenvalues are not equal to one.
If {\bf(Mon)} or an appropriate analogue for higher rank systems is satisfied then the full derived push forward reduces to the push forward with compact supports
$$Rj_*{\cal L}=j_!{\cal L}\,.$$ It follows that the Leray spectral sequence reduces to the bottom row, which is
$$E^{p,0}_2=H^p(X;j_!{\cal L})=H^p(X,D;j_*{\cal L})\,.$$
Hence
$$H^p(U;{\cal L})=H^p(X,D;j_*{\cal L})\,.$$
This is essentially the argument e.g. of \cite[Theorem 5.1]{DS}.

Now suppose that $\cal L$ is enriched, so that it belongs to the category studied in \cite{BBD} (e.g. is of geometric origin, [{\it loc.cit.} 6.2.4]). Furthermore suppose that $\cal L$ is pure of weight $w$. Then [{\it loc.cit.} 5.1.13] $Rj_*{\cal L}$ is of weight $\geq w$ and $Rj_!{\cal L}$ is of weight $\leq w$. Therefore if the condition ($\bf Mon$) is satisfied, then $Rj_*{\cal L}=j_!{\cal L}$ is pure. Under assumption that $X$ is compact the cohomology group
$H^*(X;Rj_*{\cal L})$ are pure. The admit a Hodge decomposition
$$H^k(U;{\cal L})\simeq H^k(X;Rj_*{\cal L})=\bigoplus_{p+q=k+w}H^{pq}\,.$$
The Hodge decomposition for arrangements of planes in $\P^3$ were studied in \cite{Ka}.

 If a weaker condition is satisfied: $\bf(Mon)^{*}$ of \cite[Theorem 4.1]{STV} (see also \cite{ESV}), then the local cohomology of $(R^kj_*{\cal L})_x=H^k(U\cap B_{x} {\cal L})$ is generated by logarithmic forms $df_i/f_i$, where $f_i$ is an equation of a divisor component. We would like to ask:

 \begin{question}\label{pyt}Suppose that $(R^kj_*{\cal L})_x$ is generated by logarithmic forms. Does it imply that $R^kj_*{\cal L}$ is  pure?
 \end{question}
 The sheaf $R^kj_*{\cal L}$ is a sum of locally constant sheaves supported by the intersections of divisor components. Therefore pointwise purity implies that $R^kj_*{\cal L}$ is pure as a sheaf.
 We do not give the positive answer to that Question \ref{pyt}. The passage from the analytic category to arithmetic property is delicate. We would have to compare two notions of purity: one from the world of Weil sheaves, another -- defined via differential forms, properly formulated in the sense variation of Hodge structure or mixed Hodge modules. Moreover, not every local system admits a reduction to a local system over a variety over a finite field. Work of \cite{Ti} may be a hint that something can be done for unitary local systems. When we consider local systems of geometric origin in the sense of \cite{BBD} we can freely use the notion of purity.
 Assuming purity we proceed as in the proof of Theorem \ref{t2}. The Leray spectral sequence can have nontrivial second differential, but the higher differentials vanish. Therefore
\begin{corollary} Let $U\subset X$ be an open set in a compact algebraic variety and $\cal L$ be a locally constant sheaf on $U$. Suppose that $\cal L$ is of geometric origin.
Further assume that the sheaves $R^kj_*{\cal L}$ are pure of weight $w+2k$. The skew-rows joined by the second differential of the Leray spectral sequence compute the graded pieces of the weight filtration in $H^*(U;{\cal L})$ \s
\noindent$Gr^W_{w+k+\ell} H^k(U;{\cal L})=$
 \hfill
 $$=H_{\ell}\bigl(
 0
\ot\underbrace{ H^k(X;R^0j_*{\cal L})}_0
\ot\underbrace{ H^{k-2}(X;R^1j_*{\cal L})}_{1}
\ot\underbrace{ H^{k-4}(X;R^2j_*{\cal L})}_{2}
\ot \dots
\;\bigr)\,.$$
\end{corollary}

Possibly using the methods of  \cite{Ti} one can proof degeneration of the Leray spectral sequence for the complement of hyperplane-like arrangement and a unitary local system satisfying $\bf(Mon)^*$.
We finish our discussion with the question:
\begin{question} Does in general the condition $\bf(Mon)^{*}$ of \cite[Theorem 4.1]{STV} imply degene\-ra\-tion of the Leray spectral sequence
$E^{p,q}_2=H^p(X;R^qj_*{\cal L})\;\Rightarrow\; H^{p+q}(U;{\cal L})$
degenerate?
\end{question}

\end{document}